\documentclass[a4paper]{ifacconf}

\usepackage{graphicx}      % include this line if your document contains figures
\usepackage{natbib}        % required for bibliography
\usepackage[usenames,dvipsnames]{xcolor}
\usepackage{amsmath,amssymb}
\usepackage{graphicx}
\usepackage{url}
\usepackage{float}
\usepackage{tikz}
\usepackage{pgfplots} 				% plot tikz figures
\usepackage[utf8]{inputenc}
\usepackage{dsfont}   				% to get a nice \1 
\usepackage[siunitx]{circuitikz} % Loading circuitikz with siunitx option
\usepackage[mathscr]{euscript}
\usepackage{mathtools}
\usepackage{xargs} 				% more than one optional argument
\usepackage{siunitx}
\usepackage[normalem]{ulem} % strike out text, example:  \sout{Hello}

\newcommand{\w}{\omega}
\newcommand{\1}{\mathds 1  }

\newcommand{\p}{\partial}

\newcommand{\pH}{port-Hamiltonian}

\newcommand{\G}{\mathcal{G} }
\newcommand{\V}{\mathcal{V} }

\newcommand{\E}{\mathcal{E} }
\newcommand{\N}{\mathcal{N} }

\newcommandx{\F}[2][1=1, 2=2]{\frac{#1}{#2}} % more that one optional argument

\definecolor{new}{rgb}{0.55,0,0.55}
\usepackage[prependcaption,colorinlistoftodos]{todonotes}

\DeclareMathOperator{\col}{col}

\DeclareMathOperator{\diag}{diag}

\DeclareMathOperator{\spann}{span}
\DeclareMathOperator{\rank}{rank}

\begin{document}
\begin{frontmatter}

  %%%%%%%%%%%%%%%%%%%%%%%%%%%%%%%%%%%%%%%%%%%%%%%%%% 

  \title{Stabilization of Structure-Preserving Power Networks with Market Dynamics\thanksref{footnoteinfo}}
  % \title{Style for IFAC Conferences \& Symposia: Use Title Case for   Paper Title\thanksref{footnoteinfo}} 

  %%%%%%%%%%%%%%%%%%%%%%%%%%%%%%%%%%%%%%%%%%%%%%%%%% 

  \thanks[footnoteinfo]{This work is supported by the NWO (Netherlands Organisation for Scientific Research) programme  \emph{Uncertainty Reduction in Smart Energy Systems (URSES)} under the auspices of the project \text{ENBARK}.}
  % \thanks[footnoteinfo]{Sponsor and financial support acknowledgment goes here. Paper titles should be written in uppercase and lowercase letters, not all uppercase.}

  %%%%%%%%%%%%%%%%%%%%%%%%%%%%%%%%%%%%%%%%%%%%%%%%%% 

  \author[TjerkClaudio]{Tjerk W. Stegink} 
  \author[TjerkClaudio]{Claudio De Persis} 
  \author[Third]{Arjan J. van der Schaft}

  \address[TjerkClaudio]{Engineering and Technology institute Groningen (ENTEG), University of Groningen, Nijenborgh 4, 9747 AG Groningen, \\ The Netherlands. (e-mail: \{t.w.stegink,c.de.persis\}@rug.nl)}
  \address[Third]{Johann Bernoulli Institute for Mathematics and Computer Science, University of Groningen, Nijenborgh 9, 9747 AG Groningen,\\
    The Netherlands. (e-mail: a.j.van.der.schaft@rug.nl)}

  %%%%%%%%%%%%%%%%%%%%%%%%%%%%%%%%%%%%%%%%%%%%%%%%%% 

  \begin{abstract}                % Abstract of not more than 250 words.
    This paper studies the problem of maximizing the social welfare  while stabilizing both the physical power network as well as the market dynamics. For the physical power grid a third-order structure-preserving model is considered involving both frequency and voltage dynamics. By applying the primal-dual gradient method to the social welfare problem, a distributed dynamic pricing algorithm in \pH\ form is obtained. After interconnection with the physical system  a closed-loop \pH\ system of differential-algebraic equations is obtained, whose properties are exploited to prove local asymptotic stability of the optimal point. 
  \end{abstract}

  %%%%%%%%%%%%%%%%%%%%%%%%%%%%%%%%%%%%%%%%%%%%%%%%%% 


  %%%%%%%%%%%%%%%%%%%%%%%%%%%%%%%%%%%%%%%%%%%%%%%%%% 

\end{frontmatter}
% ===============================================================================

\section{Introduction}
The future power network needs to operate reliably in the face of fluctuations resulting from distributed energy resources and the increased variability in both supply and demand. One of the feedback mechanisms that have been identified for managing this challenge is the use of real-time dynamic pricing. This feedback mechanism encourages consumers to modify their demand when it is difficult for system operator to achieve a balance between supply and demand \citep{borenstein2002dynamic}. In addition, real-time dynamic pricing allows to maximize the total \emph{social welfare} by fairly sharing utilities and costs associated with the generation and consumption of energy among the different control areas \citep{kiani_anna}. %, \cite{kiani2012hierarchicalR1}.

Many of the existing dynamic pricing algorithms focus on the economic part of optimal supply-demand matching \citep{kiani_anna,CDC2010_Stability}. However, if market mechanisms are used to determine the optimal power % generator
dispatch (with near real-time updates of the dispatch commands) dynamic coupling occurs between the market update process and the physical response of the power network dynamics \citep{alv_meng_power_coupl_market}. 

Consequently, under the assumption of market-based dispatch, it is essential to consider the stability of the coupled system incorporating both market operation and electromechanical power system dynamics simultaneously. 

% \subsubsection{Goals}
While on this subject a vast literature is already available, we focus on a more accurate and higher order model for the physical power network than conventionally used in the literature. In particular, a  structure-preserving model for the power network with a third-order order model for the synchronous generators including voltage dynamics is used. As a result, market dynamics, frequency dynamics and voltage dynamics are considered simultaneously.

% The approach taken in this paper is to model both the dynamic pricing controller as well as the physical network in a \pH\ way, emphasizing energy storage and power flow. This provides a unified framework for the modeling, analysis and control of power networks with market dynamics, with possible extensions to more refined models of the physical power network, including for example turbine dynamics.

\subsection{Literature review}
The coupling between a high-order dynamic structure-preserving power network and market dynamics has been studied before in \cite{alv_meng_power_coupl_market}. Here a fourth-order model of the synchronous generator is used in conjunction with turbine and exciter dynamics, which is coupled to a simple model describing the market dynamics. The results established in \cite{alv_meng_power_coupl_market} are  based on an eigenvalue analysis of the linearized system.  

% \subsubsection{Third-order model}
It is shown in \cite{trip2016internal} that the third-order (\emph{flux-decay}) model for the synchronous generator, as used in the present paper, admits a useful passivity property that allows for a rigorous stability analysis of the interconnection with optimal power dispatch controllers, even in the presence of time-varying demand.  
% trip turbine
In \cite{trip2015optimal} a structure-preserving power network model is considered with turbine dynamics where a similar internal-model controller is applied, which also has applications in microgrids, see \cite{de2016lyapunov}.  

% \subsubsection{Gradient-method based controllers}
Another commonly used approach to design optimal distributed controllers in power grids is the use of the primal-dual gradient algorithm \citep{arrow_gradmethod}, which has been proven useful in network flow theory \citep{feijer-paganini}. % ., see e.g.
The problem formulation varies throughout the literature on power systems, with the focus being on either the generation side \citep{AGC_ACC2014,you2014reverse}, the load side \citep{mallada2014distributed,zhao2015distributedAC,mallada2014optimal} or both \citep{zhang1achieving,zhangpapaautomatica}. 

% \paragraph{Power network as gradient dynamics}
Many of these references % A vast literature
focus on linear power system models coupled with gradient-method-based controllers \citep{AGC_ACC2014,you2014reverse,mallada2014distributed,mallada2014optimal,zhang1achieving,zhao2014design}. In these references the property that the  linear power system dynamics can be formulated as a gradient method applied to a certain optimization problem is exploited. This is commonly referred to as \emph{reverse-engineering} of the power system dynamics \citep{zhang1achieving,AGC_ACC2014,you2014reverse}. However, this approach falls short in dealing with models involving nonlinear power flows.

Nevertheless, \cite{zhangpapaautomatica,zhao2015distributedAC} show the possibility to achieve optimal power dispatch in structure-preserving power networks with nonlinear power flows using gradient-method-based controllers. On the other hand, the controllers proposed in \cite{zhangpapaautomatica} have restrictions in assigning the controller parameters and in addition require that the  topology of the physical network is a tree.

\subsection{Main contributions}
The contribution of this paper is to propose a novel energy-based approach to the problem that differs substantially from the aforementioned works. We proceed along the lines of \cite{LHMNLC,stegink2015unifying}, where a \pH\ approach to the design of gradient-method-based controllers in power networks is proposed. In those papers it is shown that both the power network as well as the controller designs admit a port-Hamiltonian representation which are then interconnected to obtain a closed-loop \pH\ system. In the present paper we extend some of these results to structure-preserving power networks. 

First it is shown that the dynamical model describing the power network as well as the market dynamics admit a \pH\ representation. Then, following \cite{LHMNLC,stegink2015unifying}, it is proven that all the trajectories of the coupled system converge to the desired synchronous solution and to optimal power dispatch.  

% \subsubsection{Physical model}
Since our approach is based on passivity and does not require to reverse-engineer the power system dynamics as a primal-dual gradient dynamics, it allows to deal with more complex nonlinear models of the power network. More specifically, the physical model for describing the power network in this paper admits nonlinear power flows and time-varying voltages, and is more accurate and reliable than the classical second-order model \citep{powsysdynwiley,kundur,sauerpai1998powersystem,firstSPM}. In addition, a distinction is made between generator nodes and loads nodes, resulting in a system of differential-algebraic equations. 

The results that are established in the present paper are valid for the case of nonlinear power flows and cyclic networks, in contrast to  \cite{zhang1achieving,AGC_ACC2014,you2014reverse,zhao2014design}, where the power flows are linearized and \cite{zhangpapaautomatica} where the physical network topology is a tree. Moreover, in the aforementioned references the voltages are assumed to be constant. 

While the third-order model for the synchronous generators has been studied before using passivity based techniques \citep{trip2016internal,stegink2015unifying}, the combination of gradient method based controllers with structure-preserving power network models is novel. In addition, the stability analysis does not rely on linearization and is based on energy functions which allow us to establish rigorous stability results. Moreover, we do not impose any restrictive condition on controller design parameters for guaranteeing asymptotic stability, contrary to \cite{zhangpapaautomatica}. 

The remainder of this paper is organized as follows. In Section \ref{sec:prel} the preliminaries are stated. Thereafter, the power system dynamics is introduced in Section \ref{sec:power-network-model} and a \pH\ representation of the system of differential-algebraic equations is given as well in this section. Then the dynamic pricing algorithm in \pH\ form is presented in Section \ref{sec:market-dynamics}. The closed-loop system is analyzed in Section \ref{sec:stab-clos-loop} and local asymptotic stability of the optimal points is proven. Finally, the conclusions and the future research directions are discussed in Section \ref{sec:conclusion}.

\section{Preliminaries}\label{sec:prel}
\subsection{Notation}
\label{sec:notations}
Given a symmetric matrix $A\in\mathbb R^{n\times n}$, we write $A>0 \ (A\ge0)$ to indicate that $A$ is a positive (semi-)definite matrix. The set of positive real numbers is  denoted by $\mathbb R_{>0}$ and likewise the set of nonnegative real numbers is denoted by $\mathbb R_{\ge 0}$. The notation $\1_n\in\mathbb R^n$ is used for the vector whose elements are equal to 1. The $n\times n$ identity matrix is denoted by $I_n$.
Given an ordered set $\mathcal I=\{i_1,i_2,\ldots,i_k\}$ and a vector $v\in\mathbb R^n, k\leq n$, then $\col_{i\in\mathcal I}\{v_i\}, \diag_{i\in\mathcal I}\{v_i\}$ denotes the $k$-column vector, respectively $k\times k$ diagonal matrix whose entries are given by $v_{i_1},v_{i_2},\ldots, v_{i_k}$. Likewise, given vectors $v_1,v_2$ then $\col(v_1,v_2):=\left[
  \begin{smallmatrix}
    v_1\\v_2
  \end{smallmatrix}
\right]$.  Let $f(x,y)$ be a differentiable function of $x\in\mathbb R^n,y\in\mathbb R^m$, then $\nabla f:=   \col( \frac{\p f}{\p x},\frac{\p f}{\p y})$ and $\nabla_xf:=\frac{\p f}{\p x}$ denotes the gradient of $f$ with respect to $x$. Given a twice-differentiable function $f:\mathbb R^n\to\mathbb R^n$ then the Hessian of $f$ evaluated at $x$ is denoted by $\nabla^2 f(x)$.  

\subsection{Differential algebraic equations}
Let us consider a system of differential-algebraic equations (DAE's) of the form 
\begin{subequations}
  \begin{align}
    \dot x&=f(x,y),\\
    0&=g(x,y),\label{eq:DAEalgcont}
  \end{align}\label{eq:DAEgeneral}
\end{subequations}
where $x\in\mathbb R^n $ and $y\in\mathbb R^m$.
\begin{defn}[\cite{de2016lyapunov}]\label{def:regular}
  Let $\mathscr D\subset \mathbb R^n\times \mathbb R^m$ be an open connected set. The algebraic equation $0=g(x,y)$ is \emph{regular} if the Jacobian of $g$ w.r.t. $y$ has constant full rank on $\mathscr D $, that is, 
  \begin{align*}
    \rank(\nabla_yg(x,y))=m\qquad \forall(x,y)\in\mathscr D.
  \end{align*}
\end{defn}
If the DAE-system \eqref{eq:DAEgeneral} is regular on $\mathscr D$ then by \cite{HillMareels1993StabilityDAE} the existence and uniqueness of solutions of \eqref{eq:DAEgeneral} in $\mathscr D$ over an interval $\mathcal I\subseteq\mathbb R_{\geq0}$ for any $(x(x_0,y_0,t),y(x_0,y_0,t))$ is guaranteed.
% \begin{rem}
%   In fact, a system that satisfies the conditions stated in Definition \ref{def:regular} is a so called \emph{Hessenberg index-1} system. Alternatively, it is often called a semi-explicit index-1 DAE. 
% \end{rem}
% \begin{rem}
%   If the DAE satisfies the conditions stated in Definition \ref{def:regular}, then by using the Implicit Function Theorem $y$ can (locally) be uniquely solved from the algebraic equation such that $y=y(x)$ (in a neighborhood) is a continuously differentiable function of $x$.  
% \end{rem}

% \subsubsection{Invariance principle for index-1 DAE systems}
By extending the usual LaSalle's invariance principle for ordinary differential equations, we obtain an invariance principle that can be used for the stability analysis of DAE's, see \cite{de2016lyapunov}.%\cite{wang2006h}.
\begin{thm}
  \label{thm:invprincDAE}
  Suppose the DAE \eqref{eq:DAEgeneral} is regular on $\mathscr D$ and $f,g$ are continuous differentiable functions. Let $(x,y)=(\bar x,\bar y)$ be an equilibrium of \eqref{eq:DAEgeneral}. Let $V(x,y):\mathscr D_V\to \mathbb R_{\geq 0}$ be a smooth positive definite function on a neighborhood $\mathscr D_V\subset \mathscr D$ of $(x,y)=(\bar x,\bar y)$, such that $\dot V(x,y)\leq 0$. Let $\mathscr S=\{(x,y)\in \mathscr D_V \ | \ \dot V=0\}$, and suppose that no solution can stay forever in $\mathscr S$, other than the trivial solution $(\bar x,\bar y)$. Then $(\bar x,\bar y)$ is locally asymptotically stable. 
\end{thm}

\section{Power network model}
\label{sec:power-network-model}
Consider a power grid consisting of $ n $ buses. The network is represented by a connected and undirected graph $ \G = (\V, \E) $. Its associated node set, $ \V = \{1, . . . , n\}=\V_g\cup \V_l $, is partitioned in the set of generator nodes $\V_g$ with cardinality $n_g$, and the set of load nodes $\V_l$ with cardinality $n_l$. The set of edges, $ \E = \{1,\ldots , m\} \subset \V \times \V $, corresponds to the set of transmission lines connecting the buses.  Each bus represents either a synchronous generator or a frequency-dependent load \citep{firstSPM}. It is assumed that the synchronous generators are governed by the \emph{flux-decay model} \citep{kundur} with a controllable mechanical power injection $P_{gi}$.  By \cite{alv_meng_power_coupl_market},  the loads are assumed to have a time-varying active power demand $P_{li}$ and a constant reactive power demand $Q_{li}=\bar Q_{li}$. As a result, the dynamics at each bus is given by  \citep{powsysdynwiley,trip2016internal}
\begin{subequations}
  \begin{align}
    \dot  \delta_i&= \w_i, && i\in \mathbb \V\\
    M_i \dot \w_i&=-P_i-A_i\w_i+P_{gi}, && i\in \mathbb \V_g\label{eq:wdyn}\\
    T_{i}\dot E_{i}&=E_{fi}-E_i-(X_{di}-X_{di}')E_i^{-1}Q_i
    , && i\in \mathbb \V_g\\
    P_i&=-A_i\w_i-P_{li}, &&i\in\mathbb \V_l,\label{eq:loaddynP}\\
    Q_i&=\bar Q_{li},    &&i\in\mathbb \V_l.\label{eq:loaddynQ}
  \end{align}\label{eq:powersysdyn}
\end{subequations}
Here the active and reactive power injections are given by %$P_{li}$ is the power demand at node $i\in\V_l$ and 
\begin{align*}
  P_i&=\sum_{j\in\N_i}B_{ij}E_iE_j\sin\delta_{ij} &&i\in\V\\
  Q_i&=B_{ii}E_i^2-\sum_{j\in\N_i}B_{ij}E_iE_j\cos\delta_{ij}, && i\in\V
\end{align*}
with $B_{ii}=\sum_{j\in\mathcal N_i}B_{ij}+\hat B_{ii}$ and where $\hat B_{ii}\geq 0$ is the negative of the shunt susceptance at node $i$.  
% which is commonly referred to as the \emph{flux-decay model}. At each load node the differential algebraic equations are given by 
% Here we use a similar notation as used in established literature on power systems \cite{powsysdynwiley}, \cite{kundur}, \cite{anderson1977}, \cite{sauerpai1998powersystem}.  
See Table \ref{tab:par3SG} for a list of symbols used in the model \eqref{eq:powersysdyn} and throughout the paper.
\begin{table}
  \centering
  \begin{tabular}[c]{ll} %\hline
    $\delta_i$     & voltage angle                                            \\ 
    $\w_i$         & frequency deviation w.r.t. nominal frequency \\%  
    $E_{i}$      & % $q$-axis
    transient internal voltage                      \\ %\hline
    $E_{fi}$       & excitation voltage                                       \\
    $P_{gi},P_{li}$       & active power generation and demand  \\%\hline
    $Q_{gi},Q_{li}$       & reactive power generation and demand        \\
    $M_i$          & moment of inertia                                        \\ 
    $\mathcal N_i$ & set of buses connected to bus $i$                        \\
    $A_i$          &  asynchronous damping constant                            \\ 
    $B_{ij}$      & negative of the susceptance of transmission line $(i,j)$ \\  
    $B_{ii}$      & negative of the self-susceptance                        \\  
    $X_{di}$       & $d$-axis synchronous reactance of generator $i$          \\
    $X_{di}'$      & $d$-axis transient reactance of generator $i$            \\
    $T_{i}$      & % $d$-axis
    open-circuit transient time constant %          \\[0.5mm]\hline
  \end{tabular}%\vspace{1mm}
  \caption{Parameters and state variables of model \eqref{eq:powersysdyn}.}  
  \label{tab:par3SG}
\end{table}
\begin{assum}\label{ass:power-network-model}
  By using the power network model \eqref{eq:powersysdyn} the following assumptions are made, where most of them are standard in a broad range of literature on power network dynamics \citep{powsysdynwiley}.
  \begin{itemize}
  \item Lines are purely inductive, i.e., the conductance is zero. This
    assumption is generally valid for the case of high
    voltage lines connecting  different control areas.
  \item The grid is operating around the synchronous frequency, say $\SI{50}{\hertz}$ or $\SI{60}{\hertz}$. 
  \item The voltages satisfy $E_i>0, i\in\V$ for all time $t\geq0$ and the reactive powers at the loads $Q_{li}^*\geq 0, i\in\V_l$ are constant. 
  \item The excitation voltage $E_{fi}$ is constant for all $i\in\V$. 
  \end{itemize}
\end{assum}
Define the angular momenta $p_i=M_i\w_i, i\in\V_g$. Based on the energy stored in the generators and the transmission network, the Hamiltonian is defined by
\begin{equation}
  \begin{aligned} H_p&=\frac1{2}\sum_{i\in\V_g}\left(M_i^{-1}p_i^2+\frac{(E_{i}-E_{fi})^2}{X_{di}-X_{di}'}\right)\\
    &+\frac1{2}\sum_{\mathclap{i\in \V}}B_{ii}E_{i}^2-\sum_{\mathclap{(i,j)\in \E}}B_{ij}E_{i}E_{j}\cos \delta_{ij}.
  \end{aligned}\label{eq:H3SG}
\end{equation}
The  algebraic constraint \eqref{eq:loaddynQ} corresponding to the reactive power of the load can then be written as 
\begin{align}\label{eq:algeqreactivepower}
  0&=-[E_l]\nabla_{E_l} H_p+\bar Q_l,
\end{align}
where $[E_l]:=\diag_{i\in\V_l}\{E_{i}\}, E_l:=\col_{i\in\V_l}\{E_i\}$.

As the system \eqref{eq:powersysdyn} admits a \emph{rotational symmetry} with respect to the voltage angles \citep{dorfler2014breaking},  it is convenient to consider a different set of coordinates. To this end,  define without loss of generality the coordinates $\varphi=\hat D^T\delta, \delta=\col_{i\in\V}\{\delta_i\}$, where  $\hat D\in\mathbb R^{n\times (n-1)}$ has full column rank and satisfies $\ker \hat D^T=\spann\{\1_{n}\}$.
\begin{rem}
  In particular, $\hat D$ could be defined as the incidence matrix of a tree graph with $n$ nodes and $n-1$ edges. In that case, $\varphi$ represents the angle differences $\delta_i-\delta_j$ along the edges of this tree graph. Another possibility is to choose node $n$ as the reference node, resulting in $\hat D^T=\left[I_{n-1} \   -\1_{n-1} \right]
  $, see also \cite{de2016lyapunov,willems1971direct}. 
\end{rem}
Using this modified set of coordinates, and defining the state variable
$x_p:=\col(\varphi,p_g,E_g,E_l)% =\col(\varphi,M \w_g,E_g,\w_l,E_l)=\tau_pz_p
$, the system \eqref{eq:powersysdyn} can be written in the form 
\begin{equation}
  \begin{aligned}
    &\begin{bmatrix}
      \dot \varphi\\
      \dot p_g\\
      \dot E_g\\
      0\\
      0
    \end{bmatrix}=
    \begin{bmatrix}
      0         & \hat D_g^T & 0      & 0 &\hat D_l^T \\        
      -\hat D_g & -A_g       & 0     & 0          & 0 \\         
      0         & 0          & -R_g   & 0          & 0 \\        
      0         & 0          & 0      &   -[E_l]         &  0   \\   
      -\hat D_l         & 0          & 0      &   0    &  -A_l
    \end{bmatrix}\nabla W_p+
    \begin{bmatrix}
      0\\P_g\\0\\\bar Q_l\\-P_l
    \end{bmatrix},\\
    &W_p(x_p,\w_l)=H_p(x_p)+U_p(\w_l), \qquad U_p(\w_l)=\frac12\w_l^T\w_l. % \qquad 
  \end{aligned}\label{eq:phphysical}
\end{equation}
Here $\hat D^T=
\begin{bmatrix}
  \hat D_g^T&\hat D_l^T
\end{bmatrix}, A_l=\diag_{i\in\V_l}\{A_i\}>0, A_g=\diag_{i\in\V_g}\{A_i\}>0, R_g=\diag_{i\in\V_g}\{\frac{X_{di}-X_{di}'}{T_i}\}>0, 
\w_l=\col_{i\in\V_l}\{\w_i\}$ and $E_g,P_g,P_l,p_g$ are defined likewise. % Observe that in the above system of DAE's the algebraic variables are the frequency deviations $\w_{li},i\in\V_l$ and the voltages $E_{i},i\in\V_l$ at the load nodes. 
The system \eqref{eq:phphysical} has external ports $(P_g,\w_g), (P_l,\w_l)$ which will be interconnected to the dynamic pricing algorithm introduced in the following section.

\begin{rem}
  The system \eqref{eq:phphysical} has a slightly different form compared to conventional \pH\ DAE-systems, see \cite{phsurvey}. In particular, $H_p$ is the Hamiltonian of the system while $U_p$ can be interpreted as an auxiliary energy function which is not used as part of the (shifted) storage function to prove passivity. However, by exploiting the special structure of the system \eqref{eq:phphysical}, the stability analysis becomes convenient as we will show in Section \ref{sec:stab-clos-loop}. 
\end{rem}

%\todoing{[R]: Change from power system model to dynamic pricing algorithm: ``The reviewer is not able to understand the proposed pricing algorithm. A comment on the transition from the DAE-system (8) to the pricing algorithm would help understanding.''}

\section{Dynamic pricing algorithm}
\label{sec:market-dynamics}
The social welfare is defined as $S(P_g,P_l):=U(P_l)-C(P_g)$, which consists of a utility function $U(P_l)$ of the power consumption $P_l$ and the  cost $C(P_g)$ associated to the power production $P_g$.  We assume that $C(P_g),U(P_l)$ are strictly convex and strictly concave functions respectively. The objective is to maximize the social welfare while achieving zero frequency deviation.  By analyzing the equilibria of \eqref{eq:powersysdyn}, it follows that a necessary condition for zero frequency deviation is $\1_{n_g}^T P_g=\1_{n_l}^T P_l$. In other words,  the total supply must match the total demand. It can be noted that $P_g,P_l$ satisfy this power balance if and only if there exists a vector $v\in\mathbb R^{m_c}$ such that
\begin{align*}-
  \underbrace{\begin{bmatrix}
      D_{cg}\\	
      D_{cl}
    \end{bmatrix}}_{D_c}v+
  \begin{bmatrix}
    P_g\\
    -P_l
  \end{bmatrix}=0
\end{align*} 
where $D_c\in\mathbb R^{n\times m_c}$ is the incidence matrix of some connected \emph{communication graph} with $m_c$ edges.  %\todoing{[R]: To what detail is the communication network considered? i.e. Does the incidence matrix account for communication latencies (=response/reaction time)?} Therefore, we will consider the following concave maximization problem:
\begin{equation}
  \begin{aligned}
    \max_{P_g,P_l,v} \ & U(P_l)-C(P_g)\\	
    \text{s.t. }	 \ &\begin{bmatrix}
      D_{cg}\\	
      D_{cl}
    \end{bmatrix}v+
    \begin{bmatrix}
      P_g\\
      -P_l
    \end{bmatrix}=0.
  \end{aligned}\label{eq:socwelfarev}
\end{equation}
The corresponding KKT optimality conditions amount to 
\begin{equation}
  \begin{aligned}
    \nabla C(\bar P_g)-\bar \lambda_g&=0\\
    -\nabla U(\bar P_l)+\bar \lambda_l&=0\\
    \begin{bmatrix}
      D_{cg}^T&D_{cl}^T
    \end{bmatrix}
    \begin{bmatrix}
      \bar \lambda_g\\
      \bar \lambda_l
    \end{bmatrix}&=0\\
    -\begin{bmatrix}
      D_{cg}\\	
      D_{cl}
    \end{bmatrix}\bar v+
    \begin{bmatrix}
      \bar P_g\\
      -\bar P_l
    \end{bmatrix}&=0
  \end{aligned}\label{eq:KKTcond}
\end{equation}
with $\lambda_g\in\mathbb R^{n_g},\lambda_l\in\mathbb R^{n_l}$. Inspired by our previous work \citep{stegink2015unifying,LHMNLC} and based on the primal-dual gradient method \citep{zhang1achieving,AGC_ACC2014,arrow_gradmethod}, the following distributed dynamic pricing algorithm is proposed:
\begin{subequations}\label{eq:dynpric}
  \begin{align}
    \tau_g\dot P_g&=-\nabla C(P_g)+\lambda_g-\w_g\\      
    \tau_l\dot P_l&=\nabla U(P_l)-\lambda_l+\w_l\\
    \tau_v\dot v&=-D_{cg}^T\lambda_g-D_{cl}^T\lambda_l \label{eq:vdyn}\\
    \tau_{\lambda_g}\dot{\lambda}_g&=D_{cg} v-P_g\label{eq:lambdag}\\
    \tau_{\lambda_l}\dot{\lambda}_l&=D_{cl} v+P_l \label{eq:lambdal}
  \end{align}
\end{subequations}
where $\tau_g,\tau_d,\tau_v,\tau_{\lambda_g},\tau_{\lambda_l}>0$ are (controller design) parameters.  %Furthermore, in the dynamic pricing algorithm \eqref{eq:dynpric} there is freedom in choosing a communication network and the associated incidence matrix $D_c$. Depending on the application, one may prefer  all-to-all communication where the underlying graph is complete, or communication networks where its associated graph is a star, line, or cycle graph.  For the present paper, it is convenient to consider a tree graph for the communication.  
Observe that the dynamics \eqref{eq:dynpric} has a clear economic interpretation \citep{kiani_anna,alv_meng_power_coupl_market}: each power producer aims at maximizing their own profit which, under the assumption of perfect competition, occurs whenever their individual marginal cost equals the local price $\lambda_{gi}-\w_{gi}$, which depends on the local frequency $\w_{gi}$ of the physical network. At the same time, each consumer maximizes its own utility but is penalized by the local price $\lambda_{li}-\w_{li}$.
\begin{rem}
  The idea to use the frequency deviation in the pricing mechanism stems from our previous work \citep{stegink2015unifying}, and helps to compensate for the power supply-demand imbalance.  Moreover, it allows for a power-preserving interconnection with the physical model \eqref{eq:powersysdyn}. 
\end{rem}
The equations \eqref{eq:vdyn},~\eqref{eq:lambdag},~\eqref{eq:lambdal} represent the distributed dynamic pricing algorithm where the quantity $v$ represents a \emph{virtual} power flow % \red{why virtual power?}
along the edges of the communication graph with incidence matrix $D_c$. We emphasize \emph{virtual}, since $v$ may not correspond to the real physical power flow as the communication graph (which can be designed as desired) may be different than the physical network topology.  By \eqref{eq:vdyn} the flow $v$ increases from areas with a lower price to areas with a relatively higher price and vice versa. Equation \eqref{eq:lambdag} shows that the local price $\lambda_i$ rises if the power outflow at node $i\in\V$ is greater than the local power supply plus power inflow of power at node $i$ and vice versa. A similar statement holds for \eqref{eq:lambdal}.  

An important observation is that the dynamic pricing algorithm  \eqref{eq:dynpric} can be written in  \pH\ form as
\begin{equation}
  \begin{aligned}
    x_c&=
    \underbrace{\begin{bmatrix}
        0 & 0  & 0   & I_{n_g}      &0\\
        0  & 0  & 0   & 0     &-I_{n_l}\\
        0  & 0  & 0   & -D_{cg}^T &-D_{cl}^T\\
        -I_{n_g}  & 0 & D_{cg} & 0   & 0\\
        0  & I_{n_l} & D_{cl} &  0  & 0\\
      \end{bmatrix}}_{J_c}\nabla H_c+
    \begin{bmatrix}
      -\w_g\\
      \w_l\\
      0\\0\\0
    \end{bmatrix}\\&+
    \nabla S(\tau_c^{-1}x_c)), \qquad
    H_c=\frac12x_c^T\tau_c^{-1}x_c 
  \end{aligned}\label{eq:phmarketpres}
\end{equation}
with 
\begin{align*}
  \underbrace{  \begin{bmatrix}
      x_g\\x_l\\x_v\\x_{\lambda_g}\\x_{\lambda_l}
    \end{bmatrix}}_{x_c}=
  \underbrace{\begin{bmatrix}
      \tau_g&0&0&0&0\\
      0&\tau_l&0&0&0\\
      0&0&\tau_v&0&0\\
      0&0&0&\tau_{\lambda_g}&0\\
      0&0&0&0&\tau_{\lambda_l}
    \end{bmatrix}}_{\tau_c}
  \underbrace{  \begin{bmatrix}
      P_g\\P_l\\v\\\lambda_g\\\lambda_l
    \end{bmatrix}}_{z_c}.
\end{align*}
Since $S$ is a concave function it satisfies the following dissipation inequality
\begin{align*}
  (z_c-\bar z_c)^T(\nabla S(z_c)-\nabla S(\bar z_c))\leq 0
\end{align*}
for all $z_c,\bar z_c\in\mathbb R^{2n+m_c}$. This property implies that the system \eqref{eq:phmarketpres} is passive with respect to its steady states, see also \cite{stegink2015unifying}. 

\section{Stability of the closed-loop system}
\label{sec:stab-clos-loop}

% \subsection{Equilibrium points}
It is observed that, by construction of the dynamic pricing algorithm, there is two-way coupling between the physical power network \eqref{eq:powersysdyn} and the market dynamics \eqref{eq:dynpric} through the ports $(P_g,\w_g), (P_l,\w_l)$. In fact, the interconnection between \eqref{eq:powersysdyn} and \eqref{eq:dynpric} is power-preserving. As a result,  the closed-loop system, obtained by combining the systems \eqref{eq:phphysical} and \eqref{eq:phmarketpres}, takes the form
\begin{equation}\label{eq:clsysDAE}
  \begin{aligned}
    \begin{bmatrix}
      \dot x_c\\
      \dot \varphi\\
      \dot p_g\\
      \dot E_g\\
      0\\
      0
    \end{bmatrix}
    &=
    \begin{bmatrix}
      J_c  & 0         & - G_1^T    & 0    & 0      & G_2^T      \\
      0    & 0         & \hat D_g^T & 0    & 0      & \hat D_l^T \\       
      G_1  & -\hat D_g & -A_g       & 0    & 0      & 0          \\        
      0    & 0         & 0          & -R_g & 0      & 0          \\       
      0    & 0         & 0          & 0    & -[E_l] & 0          \\
      -G_2 & -\hat D_l & 0          & 0    & 0      & -A_l   
    \end{bmatrix}\nabla W
    \\&+
    \begin{bmatrix}
      (\nabla S(\tau_c^{-1}x_c))^T&0&0&0&\bar Q_l^T&0
    \end{bmatrix}^T,
  \end{aligned}
\end{equation}
where $W(x,\w_l)=H_p(x_p)+U_p(\w_l)+H_c(x_c),$ with $ x=\col(x_c,x_p)$ and
\begin{align*}
  G_1=
  \begin{bmatrix}
    I_{n_g}&0&0&0&0
  \end{bmatrix}, \quad
  G_2=
  \begin{bmatrix}
    0&I_{n_l}&0&0&0
  \end{bmatrix}.
\end{align*}
Next, we examine the equilibria of the coupled system \eqref{eq:clsysDAE}. % or equivalently  \eqref{eq:powersysdyn},~\eqref{eq:dynpric}
From it follows from  \eqref{eq:vdyn} that  $\bar \lambda_g=\1_{n_g}\bar \lambda_*, \lambda_l=\1_{n_l}\lambda_*$ for some $\lambda_*\in\mathbb R$. Hence, the prices $\lambda_{gi}\in\V_g,\lambda_{li}\in\V_l$ are identical at each node. From \eqref{eq:lambdag},~\eqref{eq:lambdal} follows the power balance $\1_{n_g}^T\bar P_g=\1_{n_l}^T\bar P_l$. Finally, from \eqref{eq:wdyn} and \eqref{eq:loaddynP} we have that $\bar \w_i=0, i\in\V$. This implies that at steady state the KKT optimality conditions \eqref{eq:KKTcond} are satisfied.  Hence, the equilibrium points of the combined system \eqref{eq:powersysdyn},~\eqref{eq:dynpric} satisfy the optimality conditions of the social welfare problem \eqref{eq:socwelfarev} and, moreover, the frequency deviations are zero at steady state.

Suppose now that there exists an equilibrium $(\bar x_c,\bar x_p,\bar \w_l)$ of \eqref{eq:clsysDAE}. Then we define the shifted Hamiltonian \citep{phsurvey} by
\begin{align*}
  \bar H_p(x_p)&=H_p(x_p)-(x_p-\bar x_p)^T\nabla H_p(\bar x_p)-H_p(\bar x_p),
\end{align*}
and similarly  $\bar U_p(\w_l), \bar H_c(x_c)$ are defined. As a result, the algebraic equation \eqref{eq:algeqreactivepower} can be rewritten as 
\begin{align*}
  0&=-[E_l]\left(\nabla_{E_l}\bar H_p(x_p)+\nabla_{E_l}\bar H_a(x_p)\right)
\end{align*}
where %$\bar H_a$ is defined as
\begin{align*}
  \bar H_a(x_p)&=E_l^T\nabla_{E_l}\bar H_p(\bar x_p)-\bar Q_l^T\log E_l.
\end{align*}
with $\log(\cdot)$ being the element-wise natural logarithm. Observe that $\bar H_a$ is a convex function since, by Assumption \ref{ass:power-network-model}, $\bar Q_l\geq 0$.  % \footnote{Provided that the Hessian of $H_p$ evaluated at the equilibrium is positive definite we also allow for sufficiently small \emph{negative} reactive powers $\bar Q_l< 0$ at the loads, see also \cite{de2016lyapunov}
% .}
Since the term $\bar H_a$ only depends on $E_l$, the system \eqref{eq:clsysDAE} takes the equivalent form 
\begin{equation}\label{eq:clsysshifted}
  \begin{aligned}
    \begin{bmatrix}
      \dot x_c\\
      \dot \varphi\\
      \dot p_g\\
      \dot E_g\\
      0\\
      0
    \end{bmatrix}
    &=
    \begin{bmatrix}
      J_c    & 0         & - G_1      & 0    & 0      & G_2        \\
      0      & 0         & \hat D_g^T & 0    & 0      & \hat D_l^T \\       
      G_1^T  & -\hat D_g & -A_g       & 0    & 0      & 0          \\        
      0      & 0         & 0          & -R_g & 0      & 0          \\       
      0      & 0         & 0          & 0    & -[E_l] & 0          \\
      -G_2^T & -\hat D_l & 0          & 0    & 0      & -A_l       
    \end{bmatrix}\nabla \bar W
    \\&+
    \begin{bmatrix}
      \nabla S(\tau_c^{-1}x_c)^T-\nabla S(\tau_c^{-1}\bar x_c)^T&
      0&0&0&0&0
    \end{bmatrix}^T
  \end{aligned}
\end{equation}

where $\bar W:=\bar H+\bar U_p, \bar H:=\bar H_p+\bar H_a+\bar H_c$. 

After elimination of the algebraic variable $\w_l$, the closed-loop system \eqref{eq:clsysshifted} can in turn equivalently be  written as \eqref{eq:cleliminatedshifted}, see next page. 
\begin{figure*}[t]
  \centering
  \begin{align}\label{eq:cleliminatedshifted}
    \begin{bmatrix}
      \dot x_c                                  \\
      \dot \varphi                                \\
      \dot p_g                                  \\
      \dot E_g                                  \\
      0
    \end{bmatrix}
    & =
    \begin{bmatrix}
      J_c-G_2A_l^{-1}G_2^T   & -G_2A_l^{-1}\hat D_l         & - G_1      & 0    & 0 \\
      -\hat D_l^TA_{l}^{-1}G_2^T     & -\hat D_l^TA_l^{-1}\hat D_l         & \hat D_g^T & 0    & 0 \\       
      G_1^T & -\hat D_g & -A_g       & 0    & 0 \\        
      0     & 0         & 0          & -R_g & 0 \\       
      0     & 0         & 0          & 0    & -[E_l]       
    \end{bmatrix}\nabla \bar H    +
    \begin{bmatrix}
      \nabla S(\tau_c^{-1}x_c)-\nabla S(\tau_c^{-1}\bar x_c) \\ 
      0                                                             \\ 0          \\ 0    \\ 0
    \end{bmatrix}
  \end{align}
\end{figure*}

We are now ready to present the main convergence result.
\begin{thm}\label{thm:asymstab}
  Consider system \eqref{eq:cleliminatedshifted} and suppose that $D_c\in\mathbb R^{n\times(n-1)}$ is the incidence matrix of a tree graph. Furthermore assume that the system \eqref{eq:cleliminatedshifted} admits an equilibrium $\bar x$ satisfying $\nabla^2 \bar H(\bar x)>0$. Then $\bar x$ is locally asymptotically stable. 
\end{thm}
\begin{pf}
  The shifted Hamiltonian $\bar H$ satisfies
  \begin{equation}\label{eq:Hleq0}
    \begin{aligned}
      \dot {\bar H}&=- \w_g^TA_g \w_g- \w_l^TA_l \w_l-(\nabla_{E_g} \bar H)^TR_g\nabla_{E_g} \bar H\\&-(P_g-\bar P_g)^T(\nabla C(P_g)-\nabla C(\bar P_g))\\
      &+(P_l-\bar P_l)^T(\nabla U(P_l)-\nabla U(\bar P_l))\leq 0
    \end{aligned}
  \end{equation}
where $\w_l$ satisfies the algebraic constraint \eqref{eq:loaddynP}. Observe that $\dot {\bar H}=0$ if and only if
  $\w_g=0, \w_l=0, \nabla_{E_g}\bar H(x)=0,
  P_g=\bar P_g, P_l=\bar P_l$ since $C,U$ are strictly convex/concave functions respectively. % Let $\mathscr S=\{x\in\mathscr D \ | \ \dot {\bar H}=0\}$.
  On the largest invariant set $\Omega$ where $\dot {\bar H}=0$ we have
  \begin{alignat*}{2}
    \dot \varphi&=0            &  \lambda_g&= \bar \lambda_g\\                                               
    0   &=\hat D\nabla_\varphi \bar H &  \lambda_l&= \bar \lambda_l\\                                  
    0&=\nabla_{E_g} \bar H               &  \tau_{\lambda_g}\dot{\lambda}_g&=D_{cg} (v-\bar v) \\  
    0&=[E_l]\nabla_{E_l}\bar H     &\tau_{\lambda_l}\dot{\lambda}_l&=D_{cl} (v-\bar v)\\ 
    \tau_v\dot v&=-D_{cg}^T(\lambda_g-\bar \lambda_g)&-D_{cl}^T(\lambda_l-\bar \lambda_l&)
  \end{alignat*}
  Hence, $ \lambda_g=\bar \lambda_g, \lambda_l=\bar \lambda_l$ are constant and since $D_c,\hat D$ have full column rank it follows that $v=\bar v, \nabla_\varphi \bar H(x)=0$ respectively.  Hence, each element $x\in\Omega$ satisfies $\nabla \bar H(x)=0$. By \eqref{eq:Hleq0} and since $\nabla \bar H(\bar x)=0,\nabla^2\bar H(\bar x)>0$,  there exists a compact level set $\mathscr D$ of $\bar H$ around $\bar x$ which is forward invariant and satisfies $\mathscr D\cap \Omega=\{\bar x\}$. Hence, by Theorem \ref{thm:invprincDAE}, $\bar x$ is locally asymptotically stable.\hfill {} \qed
\end{pf}
\begin{rem}
  While the communication graph is assumed to be a tree in Theorem \ref{thm:asymstab}, the convergence result can be extended to the case of general (cyclic) connected communication graphs. However, due to space constraints, this is beyond the scope of the paper. 
\end{rem}
In \cite{de2015bregman} a sufficient condition is given which guarantees that the Hessian $\nabla^2 \bar H(\bar x)$ evaluated at the equilibrium $\bar x$ is positive definite, which is required in Theorem \ref{thm:asymstab}. Adapted for the model \eqref{eq:powersysdyn}, the following distributed condition should be verified.
\begin{prop} Suppose that $\bar x$ satisfies
  \begin{align*}
    \frac{1}{X_{di}-X_{di}'}+B_{ii}-\sum_{j\in\mathcal N_i}B_{ij}\frac{\bar E_i+\bar E_j\sin^2 \bar \delta_{ij}}{\bar E_i\cos\bar \delta_{ij}}&>0, &&i\in\V_g\\
B_{ii}-\sum_{j\in\mathcal N_i}B_{ij}\frac{\bar E_i+\bar E_j\sin^2 \bar \delta_{ij}}{\bar E_i\cos\bar \delta_{ij}}&>0, &&i\in\V_l
  \end{align*}
  with $\bar \delta_{ij}\in(-\pi/2,\pi/2), \forall (i,j)\in\E,$ and $\bar E_i>0,\forall i\in\V$. Then $\nabla^2 \bar H(\bar x)>0$.
\end{prop}
Remarkably, the above conditions are satisfied if the voltage (angle) differences and generator reactances are sufficiently  small and the shunt susceptances (at the loads) are sufficiently large.

\section{Conclusions}
\label{sec:conclusion}

In this paper an energy-based approach to the modeling and stability analysis of structure-preserving power networks with markets dynamics has been established. In particular, local convergence of the coupled system of differential-algebraic equations to the set of optimal points of the social welfare problem have been proven using a suitable Lyapunov function.

% \todoing{[T]: Perhaps mention also the extension to include constraints, following the lines of \cite{stegink2015unifying}.}

A possible extension to the established results is to consider the more complex case that the loads are not frequency dependent.  %This introduces additional challenges in the stability analysis as the closed-loop differential-algebraic system does not satisfy the regularity condition given in Definition \ref{def:regular}. As a result, the extended LaSalle's invariance principle  considered in the present paper cannot directly be applied to this case.
Another direction for future research  is to design additional controllers for the physical power network that achieve optimal reactive power sharing and/or voltage regulation. In addition, an extension can be made to include generator limits and line congestion. Finally, the influence of a possible delay in the communication of the dynamic pricing algorithm on the stability of the  closed-loop has to be investigated.  

%\todoing{[R]: Mention something about the application: ``The connection between the profound theoretical contribution of the paper and a possible application can be extended in order to help practitioners apply the theoretical results Does the proven result help in controlling power networks?''}

\bibliography{ifacconf,mybib.bib}

\end{document}